\documentclass[11pt,fleqn]{article}
\RequirePackage[T1]{fontenc}
\RequirePackage{amsthm,amsmath}
\RequirePackage[round]{natbib}
\usepackage[english]{babel}
\RequirePackage[colorlinks,citecolor=blue,urlcolor=blue]{hyperref}
\usepackage{enumerate}
\usepackage{amsmath,amsfonts,amssymb,euscript,mathrsfs}
\usepackage{fullpage}
\usepackage{amsmath}
\usepackage{dsfont}
\makeatletter
\def\captionof#1#2{{\def\@captype{#1}#2}}
\makeatother
\def\1{\mbox{\bf 1}}
\def\R{\mathbb{R}}
\def\B{\mathbb{B}}
\def\N{\mathbb{N}}
\def\P{\mathbb{P}}
\def\E{\mathbb{E}}
\def\L{\mathbb{L}}

\def\R{\mathbb{R}}

\def\Z{\mathbb{Z}}

\newtheorem{theo}{Theorem}
\newtheorem{lem}{Lemma}

\newtheorem{Def/Prop}{Definition-Proposition}

\newcounter{exos}
\renewcommand\theexos{\arabic{exos}}

\newcounter{prob}
\renewcommand\theprob{\arabic{prob}}

\begin{document}
\author{ Zinsou Max Debaly \footnote{CREST-ENSAI, UMR CNRS 9194, Campus de Ker-Lann, rue Blaise Pascal, BP 37203, 35172 Bruz cedex, France.}\and
Lionel Truquet \footnote{CREST-ENSAI, UMR CNRS 9194, Campus de Ker-Lann, rue Blaise Pascal, BP 37203, 35172 Bruz cedex, France. {\it Email: lionel.truquet@ensai.fr}}}

\title{Stationarity and Moment Properties of some Multivariate Count Autoregressions}
\date{}
\maketitle

\begin{abstract}
\noindent
We study stationarity and moments properties of some count time series models  
from contraction and stability properties of iterated random maps. Both univariate and multivariate processes are considered, including the recent multivariate count time series models introduced recently by \cite{DFT}. We improve many existing results by providing optimal 
stationarity conditions or conditions ensuring existence of some exponential moments. 
\end{abstract}
\vspace*{1.0cm}

\footnoterule

\section{Introduction}
Count time series have been widely studied in literature. 
The recent textbook of \citet{Weiss} presents many possible stochastic models for such series which offer challenging theoretical problems due to their discrete nature. Models called observation-driven, following the classification of \citet{Cox1981}, are one of the most popular. In such models, the conditional distribution at time $t$, $X_t\vert X_{t-1},X_{t-2},\ldots$ depends on an unobserved random parameter $\lambda_t$ defined by some recursions 
$\lambda_t=G_{\theta}\left(\lambda_{t-1},Y_{t-1},\ldots,\lambda_{t-q},Y_{t-q}\right)$, with a measurable function $G_{\theta}$ depending on an unknown parameter $\theta\in\R^d$. 
A typical example are INGARCH processes, defined from a Poisson distribution and a linear function $G_{\theta}$. See for instance \citet{Ferland} and \citet{Fok} for basic properties of such processes.
Though theoretical properties of such models have been extensively studied in the literature, it is difficult to find a general approach to study their stationarity properties. One difficulty arises using the standard Markov chain techniques. Though $(Y_t,\lambda_t)$ form a Markov process, the usual irreducibility assumption
is not satisfied due to the second component. This led some authors to use perturbation methods as in \citet{Fok}, \citet{Fokianos2011b} or \citet{EJS}, in order to approximate the dynamic of the non-irreducible Markov chain by an irreducible one. Another widely used approach is based on contraction techniques for autoregressive processes. However, though the main idea is to use contraction on average methods, the proposed frameworks are never unified. 
For instance, \citet{Ferland} use approximation of INGARCH processes with INAR processes, \citet{doukhan2012weak} use representation of $Y_t$ as an infinite memory process,  \citet{davis} use contraction for iterated random maps developed by \citet{Wu2}. 
With the same type of approach, a few contributions also consider multivariate time series of counts. See in particular \citet{Latour} for multivariate INAR processes or the recent contribution of \citet{DFT} for extension of INGARCH processes to the multidimensional case.  
Our motivation with this paper is to use contraction and stability properties of some iterated random maps to define a unified framework for studying univariate and multivariate time series of counts, when contraction techniques are relevant to study such models. At this point, we precise that we will not introduce a new sophisticated method but we simply reformulate some basic ideas to get a general result which is simple to apply and is quite well adapted to the standard INGARCH and INAR process as well as to their multivariate extensions.
With our formulation, we recover some well-known results but also new ones. In particular we improve several results given in \citet{DFT} and we study existence of exponential moments for INGARCH processes and its multivariate extension. Existence of exponential moments are particularly important for studying consistency of regularized estimators in multivariate time series models. See for instance \citet{LASSO} for some assumptions ensuring consistency of the LASSO estimator in general time-dependent regression models. 
       
The paper is organized as follows. In Section \ref{GEN}, we recall a result of \citet{Wu2} about random maps contracting on average and we give an extension of this result to higher-order autoregressive processes that will be particularly useful in the multivariate case. We also provide a simple result which ensures existence of some moments. Applications to some multivariate count autoregressions are considered in Section \ref{MOD}. Finally, two technical lemmas are given in an appendix section. 

\section{General result}\label{GEN}
In this section, we recall a standard result given in \citet{Wu2} about the convergence of the backward iterations of random maps. We next provide and extension to $q$-order multivariate autoregressive models which will be particularly useful for studying some multivariate counts autoregressions. We also give some results which guaranty existence of some moments for the stationary solution.
\subsection{Iterated random maps system}
Let $E$ be a subset of $\R^k$, $k \in \N^{*}$. We equipped $E$ with the Borel algebra $\mathcal{E}.$  The $1-$order autoregressive system on $(E, \mathcal{E})$ is defined as a sequence of random variables $(X_t)_{t \in \Z}$ taking values in $(E, \mathcal{E})$ such that : 
\begin{equation}\label{musc}
 \forall t \in \mathbb{Z}, X_{t+1} = F_{\epsilon_{t+1}}(X_{t}), 
    \end{equation}
where $(\epsilon_t)_{t\in \Z}$ is a sequence of random variables taking values in a second measurable space $(G, \mathcal{G})$ and independent and identically distributed (i.i.d.). The application $(x,s) \mapsto F_s(x)$  is assumed to be measurable, as an application from $E\times G$ to $E$. The equation (\ref{musc}) assumes that there exists a background process $\left\{f_t:=F_{\epsilon_t}(\cdot): t\in\Z\right\}$ that transforms the values of sequence $(X_t)_{t \in \Z}$ at time $t$ to the next one.  The following theorem, which is given in \cite{Wu2}, provides sufficient conditions for existence of a stationary solution for (\ref{musc}). 
For notational convenience, we set for two integers $s<t$, $f_s^t(x)=f_t\circ f_{t-1}\circ\cdots\circ f_s(x)$.
The following result can be found in \citet{Wu2}, Theorem $2$.

\begin{theo}
\label{th::main}
We assume that  there exists a norm $\vert\cdot\vert$ on $E$, $C>0$ and  $\kappa \in (0,1)$ and an integer $m$ such that the following assumptions are satisfied.
\begin{enumerate}[(H1)]
    \item For $x\in E$, the random variable $f_0(x)$ has a finite first-order moment:
    $$
    \mathbb{E}(\vert f_0(x)\vert)< \infty.
    $$
    \item For $(x,y)\in E^2$, we have:
    $$
    \mathbb{E}(\vert f_0(x)-f_0(y)\vert) \leq C \vert x-y\vert.
    $$
    \item For $(x,y)\in E^2$, 
				$$\mathbb{E}(\vert f_1^m(x)-f_1^m(y)\vert) \leq \kappa \vert x-y\vert.$$
\end{enumerate}
Then, for any $x \in E$ and $t\in\Z$,  the backward iterations $f_{t-n}^t(x)$ 
converge, as $n\rightarrow\infty$, almost surely and in $\L^1$ to a random variable $X_t$ not depending on $x$.
The process $(X_t)_{t\in\Z}$ is stationary and ergodic, solution of (\ref{musc}). Moreover, $\P_{X_0}$ is the unique invariant probability measure of a Markov chain defined by (\ref{musc}). 
\end{theo}

We now provide a useful lemma which guaranty existence of some moments without requiring contraction properties. 
The following result is a standard application of the concept of drift function, usually denoted by $V$, widely used for studying stability properties of Markov chains.
See \citet{MT}. Two important cases in what will follow are $V(x)=\vert x\vert^r$ for some $r>1$ and $V(x)=\exp(\delta x)$ for $\delta>0$.
\begin{lem}
\label{lem::moment}
Suppose that {\bf (H1)-(H3)} hold true and that $V:E\rightarrow (0,\infty)$ is a drift function, i.e. there exists $b>0$, $\overline{\kappa}\in (0,1)$ and
an integer $\overline{m}\geq 1$ such that for all $x\in E$,
$$\E\left[V\left(f_1^{\overline{m}}(x)\right)\right]\leq \overline{\kappa}V(x)+b.$$
If $V$ is continuous, then $\E\left[V(X_0)\right]<\infty$.  
\end{lem}

\paragraph{Proof of Lemma \ref{lem::moment}}
Using the almost sure convergence of $\left(f_{t-n\overline{m}}^t(x)\right)_{n\geq 1}$ and Fatou's lemma, we have
$$\mathbb{E}\left[V(X_{-1})\right]\leq \liminf_{n\rightarrow \infty}\mathbb{E}\left[V\left(f_{-n\overline{m}}^{-1}(x)\right)\right].$$
Using the assumption on $V$, the decomposition $f_{-n\overline{m}}^{-1}(x)=f_{-(n-1)\overline{m}}^{-1}\circ f_{-nm}^{-(n-1)\overline{m}-1}(x)$ and the fact that the $f_t'$s are i.i.d., we have 
$$\mathbb{E}\left[V\left(f_{-n\overline{m}}^{-1}(x)\right)\right]\leq \overline{\kappa}\E\left[V\left(f_{-(n-1)\overline{m}}^{-1}(x)\right)\right]+b.$$
Iterating the previous bound, we get 
$$\mathbb{E}\left[V\left(f_{-n\overline{m}}^{-1}(x)\right)\right]\leq \sum_{j=0}^{n-1} \overline{\kappa}^j b +\overline{\kappa}^nV(x).$$
Since the right hand side of the last inequality is bounded in $n$, we get integrability of $V(X_0)$.$\square$

\subsection{Result for higher-order autoregressive processes}
We generalize the previous results to some higher-order random iterations. 
Our results, which will be useful for non linear multivariate autoregressions, differ from the standard strong contraction assumptions available in the literature such as in \citet{Duflo} or \citet{DW} among others. In particular, we use a vectorial formulation which allows contraction properties for iterated nonlinear random maps.  
 
In this section, we still consider 
a sequence $(\epsilon_t)_{t\in\Z}$ \textit{i.i.d} of random variables taking values in a measurable space $G$, a subset $E$ of $\R^k$ endowed with a norm $\vert\cdot\vert$. Our aim is to study existence of stationary solutions for 
the following recursive equations:
\begin{equation}
\label{eq::recur}
X_t=F\left(X_{t-1},\ldots,X_{t-q},\epsilon_t\right),\quad t\in\Z,
\end{equation}
where $F:E^q \times G\rightarrow E$ is a measurable function. 
We first introduce additional notations. For any positive integer $d$, we denote by $\mathcal{M}_d$ the set of square matrices with real coefficients and $d$ rows and if $A\in\mathcal{M}_d$, $\rho(A)$ the spectral radius of the matrix $A$. Moreover,
for $x\in \R^d$ and $r\in\R_+$, the vector $\left(\vert x_1\vert^r,\ldots,\vert x_q\vert^r\right)'$ will be denoted by $\vert x\vert_{vec}^r$. Finally, we introduce a partial order relation $\preceq$ on $\R^d$ and such that $x\preceq x'$ means $x_i\leq x_i'$ for $i=1,\ldots,d$.

The following assumptions will be needed.
\begin{description}
\item [A1]
For any $y\in E^q$, $\E\left[\left\vert F\left(y,\epsilon_0\right)\right\vert\right]<\infty$.
\item [A2] 
There exists some matrices $A_1,\ldots,A_q\in \mathcal{M}_p$ with nonnegative elements, satisfying $\rho\left(A_1+\cdots+A_q\right)<1$ and  such that for $y,y'\in E^q$,
$$\E\left[\left\vert F(y,\epsilon_1)-F(y',\epsilon_1)\right\vert_{vec}\right]\preceq \sum_{i=1}^q A_i\left\vert y_i-y'_i\right\vert_{vec}.$$
\item [A3] For an integer $\overline{k}>1$, there exist a vector $\phi:=\left(\phi_1,\ldots,\phi_{\overline{k}}\right)$ of continuous functions from $E\rightarrow \R_+$, a real number $r\geq 1$, some matrices $D_1,\ldots,D_q\in \mathcal{M}_{\overline{k}}$ with nonnegative elements such that $\rho\left(D_1+\cdots+D_q\right)<1$ and $c\in\R_+^p$ such that for $y\in E^q$, 
$$
\left\|\phi\left(F(y, \epsilon_1)\right) \right\|_{r, vec} \preceq c + \sum_{i=1}^q D_i \phi(y_i),
$$
where for  a random vector $Z=\left(Z_1,\ldots,Z_{\overline{k}}\right)$, $\left\| Z \right\|_{r, vec} := (\E^{1/r}\left[\vert Z_1\vert^r\right], \ldots,  \E^{1/r}\left[\vert Z_{\overline{k}}\vert^r\right])'$.
\end{description}

In what follows, for any positive integer $k$ and $x=(x_1,\ldots,x_k)\in\R^k$, we set $\vert x\vert_1=\sum_{i=1}^k\vert x_i\vert$. $\vert\cdot\vert_1$ is then the $\ell_1-$norm on $\R^k$.
\begin{theo}
\label{th::autoreg}
Assume that Assumptions {\bf A1-A2} hold true.
\begin{enumerate}
    \item There then exists a unique stationary and
    non-anticipative process $(X_t)_{t\in\Z}$ solution of \eqref{eq::recur} such that $\E\left[\vert X_t\vert\right]<\infty$. 
    \item If in addition, {\bf A3} holds true, then $\E\left[\vert \phi(X_0)\vert_1^r\right]<\infty.$
		\end{enumerate}
\end{theo}
By non-anticipative, we mean that $X_t$ is measurable with respect to $\sigma\left(\epsilon_s:s\leq t\right)$.

\paragraph{Proof of Theorem \ref{th::autoreg}}
\begin{enumerate}
\item
Define the following random map 
$$f_t(u_1,\ldots,u_q)=\left(F\left(u_1,\ldots,u_q,\epsilon_t\right)',u'_1,\ldots,u'_{q-1}\right)'$$
and the sigma fields $\mathcal{F}_t=\sigma\left(\epsilon_s:s\leq t\right)$, $t\in \Z$.
We first note that a process $(U_t)_{t\in\Z}$ satisfies (\ref{eq::recur}) if and only if the process
$(X_t)_{t\in\Z}$ defined by $X_t=\left(U_t,\ldots,U_{t-q+1}\right)'$ satisfies the recursions $X_t=f_t(X_{t-1})$, $t\in\Z$.
It then only remains to study existence of stationary solutions for the recursions defined by the random functions $f_t$, $t\in\Z$. 
We set $x=(u_1,\ldots,u_q)\in E^q$ and for $1\leq t\leq q$, $U_t(x)=u_{q-t+1}$. Next for $t\geq q+1$, we define $U_t(x)$ recursively by 
$$U_t(x)=F\left(U_{t-1}(x),\ldots,U_{t-q}(x),\varepsilon_t\right).$$ 
We then have for $t\geq q+1$,
$$\left(U_t(x),\ldots,U_{t-q+1}(x)\right)=f_{q+1}^t(x).$$
Using our assumptions, we have for $t\geq q+1$, 
$$\E\left[\vert U_t(x)-U_t(x')\vert_{vec}\vert \mathcal{F}_{t-1}\right]\preceq \sum_{i=1}^q A_i\left\vert U_{t-i}(x)-U_{t-i}(x')\right\vert_{vec}.$$
Setting $w_t=\E\left[U_t(x)-U_t(x')\right]$ for $t\geq 1$, we have $w_t\preceq \sum_{i=1}^q A_i w_{t-j}$ for $t\geq q+1$. From Lemma \ref{calc} given in the Appendix,  there exist constants $C>0$ and $\rho\in (0,1)$ such that 
$$\vert w_t\vert_1\leq C\rho^t\vert x-x'\vert_1.$$
Next, we set 
$$V_t(x)=\left(U_t(x)',\ldots,U_{t-q+1}(x)'\right),\quad t\geq q+1.$$
We have $V_t(x)=f_{q+1}^t(x)$. Checking {\bf H1-H2} is straightforward whereas Assumption {\bf H3} is satisfied with the $\ell_1-$norm and for $m$ large enough to have $qC'\rho^{m+q}<1$. Theorem \ref{th::main} ensures existence of stationary solution.

Next, let $(X_t)_{t\in\Z}$ and $(Y_t)_{t\in\Z}$ be two non-anticipative stationary solutions of (\ref{eq::recur}), both having a finite first moment.
We have for any $t\in\Z$,
$$\E\left[\left\vert X_t-Y_t\right\vert_{vec}\right]\preceq \sum_{j=1}^q A_j \E\left[\left\vert X_{t-j}-Y_{t-j}\right\vert_{vec}\right].$$
Using Lemma \ref{calc}, we deduce that $\E\left[\left\vert X_t-Y_t\right\vert_{vec}\right]=0$ and then $X_t=Y_t$ a.s. This shows the uniqueness.
\item
Using the notations of the previous point, we have 
$$\Vert \phi\left(U_t(x)\right)\Vert_{r,vec}\preceq c+\sum_{i=1}^q D_i \Vert \phi\left(U_{t-i}(x)\right)\Vert_{r,vec}.$$
Using Lemma \ref{calc} (2.) and the triangular inequality, we get 
$$\E^{1/r}\left[\vert \phi\left(U_t(x)\right)\vert^r_1\right]\leq C\overline{\rho}^t\sum_{i=1}^q\vert \phi(u_i)\vert_1+D,$$
for some constants $C,D>0$ and $\overline{\rho}\in (0,1)$ and only depending on $c,D_1,\ldots,D_q$.
Remembering that $f_{q+1}^t(x)=\left(U_t(x)',\ldots,U_{t-q+1}(x)'\right)'$ and
setting $V(x)=\left(\sum_{i=1}^q\vert\phi(u_i)\vert_1\right)^r$, the assumptions of Lemma \ref{lem::moment} are satisfied for large values of $t$. The result then follows from Lemma \ref{lem::moment}. $\square$

\end{enumerate}

\section{Multivariate count autoregressions}\label{MOD}

\subsection{Notations}
We first introduce additional notations. Let $\vert\cdot\vert_1$, $\vert\cdot\vert_2$ and $\vert\cdot\vert_{\infty}$ be the three norms on $\R^p$ defined by 
$$\vert x\vert_1=\sum_{i=1}^p\vert x_i\vert,\quad \vert x\vert_2=\sqrt{\sum_{i=1}^px_i^2},\quad \vert x\vert_{\infty}=\max_{1\leq i\leq p}\vert x_i\vert.$$
We still denote by $\vert\cdot\vert_{j}$, $j=1,2,\infty$, the associated operator norms on the space $\mathcal{M}_p$ of square matrices $p\times p$ and with real-valued coefficients. More precisely,  
$$\vert A\vert_j=\sup_{x\neq 0}\frac{\vert A x\vert_j}{\vert x\vert_j},\quad A\in\mathcal{M}_p.$$
We remind that 
$$\vert A\vert_1=\max_{1\leq j\leq p}\sum_{i=1}^p\vert A(i,j)\vert,\quad \vert A\vert_{\infty}=\max_{1\leq i\leq p}\sum_{j=1}^p\vert A(i,j)\vert.$$

Finally, for a column vector $x$ with $p$ coordinates, we set $\vert x\vert_{vec}=\left(\vert x_1\vert,\ldots,\vert x_p\vert\right)'$ and is $A$ is a matrix, $\vert A\vert_{vec}=\left(\vert A(i,j)\vert\right)_{i,j}$.

\subsection{Multivariate GINAR$(q)$ process}

This model has been studied by \citet{Latour} and using Theorem \ref{th::autoreg}, we will recover many results but also get additional moment properties. 
The model writes 
\begin{equation}\label{INAR}
X_t=\sum_{j=1}^q A_{t,j}\circ X_{t-j}+U_t,\quad t\in\Z,
\end{equation}
where for $x\in\N^p$,
$$A_{t,j}\circ x=\left(\sum_{\ell=1}^pA_{t,j}(i,\ell)\circ x_{\ell}\right)_{1\leq i\leq p}$$
and for $y\in \N$,
$$A_{t,j}(i,\ell)\circ y=\sum_{s=1}^y Y^{t,j,i,\ell}_s.$$
The latter operator $\circ$ is called the thinning operator.
We assume that $(U_t)_{t\in\Z}$ is a sequence of i.i.d. integrable random vectors in $\N^p$ and independent from the family
$$\left\{Y^{t,j,i,\ell}_s:(s,t,j,i,\ell)\in \Z^2\times \{1,\ldots,q\}\times \{1,\ldots,p\}^2\right\}$$
which is itself composed of independent integrable and integer-valued random variables and such that for $(s,t,s',t',j,i,\ell)\in \Z^4\times \{1,\ldots,q\}\times \{1,\ldots,p\}^2$, $Y^{t,j,i,\ell}_s$ et $Y^{t',j,i,\ell}_{s'}$ have the same distribution with mean $A_j(i,\ell)$.
When $q=1$, this process coincides with a Galton-Watson process with immigration. 
Note that for a process $(X_t)_{t\in\Z}$ defined by (\ref{INAR}), we have 
$$\E\left[X_t\vert X_{t-1},\ldots,X_{t-q}\right]=\sum_{j=1}^q A_j X_{t-j}+\E[U_0].$$

We will use the two following assumptions.

\begin{description}
\item[G1] The spectral radius of the matrix $A_1+\cdots+A_q$ is less than $1$.
\item[G2] There exists $r>1$ such that for all $(j,i,\ell)\in\{1,\ldots,q\}\times\{1,\ldots,p\}^2$, $Y_0^{0,j,i,\ell}$ and $U_0$ have a moment of order $r$.
\end{description}

\begin{theo}\label{bbin}
Assume that Assumption {\bf G1} holds true. There then exists a unique stationary, non-anticipative and integrable solution to the recursions (\ref{INAR}).
If in addition, Assumption {\bf G2} is valid, we have $\E\vert X_0\vert^r_1<\infty$.
\end{theo}

\paragraph{Note.} We obtain the same result as \citet{Latour} for the existence of a square integrable stationary solution for the recursions (\ref{INAR}). 
However, using our formalism, we avoid lengthy computations to check such results.
We also provide conditions for existence of a moment of arbitrary order $r>1$, a problem not investigated in \citet{Latour}.

\paragraph{Proof of Theorem \ref{bbin}}
The noise at time $t$, denoted by $\epsilon_t$, is a vector with components $U_t$ and random sequences $Y^{t,j,i,\ell}_{\cdot}$ for $1\leq j\leq q$, $1\leq i,\ell\leq p$. Define $F(x_1,\ldots,x_q,\epsilon_t)=\sum_{j=1}^q A_j\circ x_j+\epsilon_t$.
Due to the properties of the thinning operator, we have for $x_1,\ldots,x_q,y_1,\ldots,y_q\in\N^p$,
$$\E\left\vert F(x_1,\ldots,x_q,\epsilon_0)-F(y_1,\ldots,y_q,\epsilon_0)\right\vert_{vec}\preceq \sum_{j=1}^q A_j \vert x_j-y_j\vert_{vec}.$$ 
This shows Assumption {\bf A2}. Assumption {\bf A1} is automatically satisfied. The first part of the theorem follows from Theorem \ref{th::autoreg} (1.)  

Finally, we check {\bf A3} when $\phi$ is the identity function. We decompose
$$F(x_1,\ldots,x_q,\epsilon_0)=\sum_{j=1}^q A_j x_j+S(x_1,\ldots,x_q)+U_0,$$
where $S(x_1,\ldots,x_q)=\sum_{j=1}^q A_j\circ x_j-\sum_{j=1}^q A_j x_j$ is a vector of sums of independent random variables. Using Burkh\"older's inequality 
we have for $y\in\N$, $\Vert A^{0,j,i,\ell}_0\circ y-A_j(i,\ell)y\Vert_r\leq C y^{1/\max(r,2)}$ where $C>0$ depends on $r$ and $\L^{r}-$norm of the counting sequences.
One can then take the same constant $C$ for all the counting sequences.
We denote by $H$ the matrix $p\times p$ with all components equal to $1$.
For any $\varepsilon>0$, there exists $b_{\varepsilon,r}>0$ only depending on $\varepsilon,r$ and such that for $y\in\N$, $y^{1/\max(r,2)}\leq \varepsilon y+b_{\varepsilon,r}$. We choose $\varepsilon>0$ such that $\rho\left(A_1+\cdots+A_q+C\varepsilon q H\right)<1$.
We then obtain
\begin{eqnarray*}
\Vert F(x_1,\ldots,x_q,\epsilon_0)\Vert_{r,vec}&\preceq& \sum_{j=1}^q A_j x_j+\Vert S(x_1,\ldots,x_q)\Vert_{r,vec}+\Vert U_0\Vert_{r,vec}\\
&\preceq& \sum_{j=1}^q \left(A_j+C\varepsilon q H\right)x_j+qpCb_{\varepsilon,r}\mathds{1}+\Vert U_0\Vert_{r,vec}.
\end{eqnarray*}
Setting $D_j=A_j+C\varepsilon q H$ for $j=1,\ldots,q$, the result follows from Theorem \ref{th::autoreg} (2.)$\square$

\subsection{Multivariate INGARCH model with linear intensity}
We consider the multivariate count process $(Y_t)_{t\in\Z}$ introduced recently by \citet{DFT}. It is defined by 
\begin{equation}\label{multlin}
Y_t=N^{(t)}_{\lambda_t},\quad \lambda_t=d+\sum_{i=1}^q A_i \lambda_{t-i}+\sum_{i=1}^q B_i Y_{t-i},
\end{equation}
where the $A_i'$s and the $B_i'$s are $p\times p$ matrices of nonnegative elements, $d$ is a vector of $\R_+^p$
and $\left(N^{(t)}\right)_{t\in\Z}$ is a sequence of i.i.d. $p-$dimensional count processes such that 
$N_1^{(t)},\ldots,N_p^{(t)}$ are homogenous Poisson processes with intensity $1$. 
For $\lambda\in\R_+^p$, the $p-$dimensional random vector $N^{(t)}_{\lambda}$ is equal to $\left(N^{(t)}_{1,\lambda_1},\ldots,N^{(t)}_{p,\lambda_p}\right)'$.
Let us note that when $p=1$, the process coincides with the INGARCH model developed by \citet{Ferland} and \citet{Fok} and for which the conditional distribution of $Y_t$ given past values is a Poisson distribution with random intensity $\lambda_t$. 
For $p>1$, a particular case of (\ref{multlin}) is obtained when $N^{(t)}$ is a vector of independent Poisson processes. 
\citet{DFT} provided a more general approach using copula. While all the coordinates of $N^{(t)}$ are still Poisson processes with intensity $1$, they can have a quite general dependence structure. However, all the results given in \citet{DFT} about existence of stationary solutions and their marginal moments are independent from this dependence structure. This will be also the case for the results given in the present paper.

The following result provides a necessary and sufficient condition for existence of a stationary solution for (\ref{multlin}). 
\begin{theo}
\label{pp::lm}
Assume that $\rho(\sum_{i=1}^q (A_i + B_i)) < 1$. 
\begin{enumerate}
\item
There then exists a unique non anticipative, stationary and integrable solution $(Y_t)_{t\in\Z}$ for (\ref{multlin}).
\item
If (\ref{multlin}) admits a stationary solution and all the components of $d$ are positive, then $\rho\left(\sum_{i=1}^q(A_i+B_i)\right)<1$.
\item
For any $r>1$, we have $\E\left[\vert Y_t\vert_1^r\right]<\infty$.
\item
In contrast, assume that $\sum_{i=1}^q\left(\vert A_i\vert_1+\vert B_i\vert_1\right)<1$ or $\left\vert \sum_{i=1}^q(A_i+B_i)\right\vert_{\infty}<1$. There then exists $\delta>0$ such that 
$\E\left[\exp\left(\delta\vert Y_0\vert_1\right)\right]<\infty$ and $\E\left[\exp\left(\delta \vert \lambda_0\vert_1\right)\right]<\infty$. 
\end{enumerate}
\end{theo}

\paragraph{Notes}
\begin{enumerate}
\item
When $q=1$, \citet{DFT} provided various sufficient conditions for existence of a stationary solution for (\ref{multlin}). 
In contrast, Theorem \ref{pp::lm} provides an optimal condition for stationarity. In particular, the condition on the spectral radius $\rho\left(\sum_{i=1}^q(A_i+B_i)\right)<1$ is implied by any contraction condition of the form $\left\vert \sum_{i=1}^q(A_i+B_i)\right\vert<1$, where $\vert\cdot\vert$ is a matrix norm. 
We remind that for any operator norm $\vert\cdot\vert$ and any matrix $A$ of size $p\times p$, we have $\vert A\vert\leq \rho(A)$. 
Contraction conditions derived from the matrix norms $\vert\cdot\vert_1$ or $\vert\cdot\vert_2$ are used in \citet{DFT}. 
To enlighten the difference, consider the case $q=1$, $A_1=0$ and $B_1=\begin{pmatrix} \alpha&\beta\\0&\alpha\end{pmatrix}$. Conditions $\vert B_1\vert_1<1$ or $\vert B_1\vert_{\infty}$ means $\alpha +\beta<1$, while $\vert B_1\vert_2\geq \sqrt{\alpha^2+\beta^2}$. In contrast $\rho(B_1)=\alpha$ and the condition $\alpha<1$ is a substantial improvement of the restrictions obtained from contractions with respect to the previous norms. 
\item
The last point of Theorem \ref{pp::lm} provides conditions for existence of exponential moments. Assume that $p=1$.  
In the univariate case, the various contraction conditions are equivalent to $\sum_{i=1}^q(A_i+B_i)<1$ which is an optimal condition for existence of a stationary and integrable solution. Under this assumption, existence of polynomials moments have been widely discussed in the literature. See for instance \citet{Ferland} or \citet{Fok}. Theorem \ref{pp::lm} provides a stronger result by showing that this stationarity condition is sufficient for existence of exponential moments. 
\item
 A challenging question in the case $p>1$ is the following: do we still have exponential moments only using the stationarity condition $\rho\left(\sum_{i=1}^q(A_i+B_i)\right)<1$ ? While Theorem \ref{pp::lm} (3.) shows that the latter condition entails existence of moments of any order, finiteness of exponential moments is unclear. 
Note that we have an affirmative answer in the univariate case.   
\end{enumerate} 

\paragraph{Proof of Theorem \ref{pp::lm}}
We set $E=\N^p\times \R_+^p$.
\begin{enumerate}
\item
We first note that any solution $X_t = (Y_t, \lambda_t)$ of the problem satisfies the recursions
$$X_t = F(X_{t-1}, \ldots, X_{t-q}, N^{(t)}),$$
where $\forall j = 1, \ldots, q, x_j = (y_j, s_j)\in (\N^p \times \R_+^p)$,
$$F\left(x_1,\ldots,x_q,N^{(t)}\right)=\left(N^{(t)}_{f(x_1, \ldots, x_q)},f(x_1, \ldots, x_q) \right)',$$
$$f(x_1, \ldots, x_q) = d+\sum_{j=1}^q A_j s_{j}+\sum_{j=1}^q B_j y_{j}.$$

For $x \in E^q$, we have the equality 
      $$
      \E\left[\left\vert F\left(x,N^{(1)}\right)\right\vert_1\right] = 2 \mathds{1}'\left(d+\sum_{j=1}^q A_j s_{j}+\sum_{j=1}^q B_j y_{j}\right) < \infty.
      $$
 
 Moreover, for $x,x' \in E^q$ with $x = (x_1, \ldots, x_q), x' = (x_1', \ldots, x_q')~\forall j = 1, \ldots, q, x_j = (y_j, s_j), x_j = (y_j', s_j')$, 
 $$\E\left[\left\vert F(x,N^{(1)})-F(x',N^{(1)})\right\vert_{vec}\right]\preceq \sum_{j=1}^q \begin{pmatrix}
B_j  & A_j \\
B_j & A_j 
 \end{pmatrix} \left\vert x_j-x'_j\right\vert_{vec}.$$
One can notice  that the matrices  $ \Gamma = \sum_{j=1}^q \begin{pmatrix}
B_j & A_j \\
B_j & A_j 
 \end{pmatrix}$ and $\sum_{j=1}^q 
(A_j + B_j)
$ have the same spectral radius. The result then follows from Theorem \ref{th::autoreg}.

\item
If $(Y_t)_{t\in\Z}$ is a stationary and integrable solution of (\ref{multlin}), we have 
$$\E(Y_t)=\E(\lambda_t)=d+\sum_{j=1}^q B_j \E(Y_{t-j})+\sum_{j=1}^q A_j \E(\lambda_{t-j}).$$
 Setting $m=\E(Y_t)$, we have $m=d+Em$ with $E=\sum_{j=1}^q(A_j+B_j)$. 
We then obtain $m=d+Ed+\cdots+E^{n-1}d+E^n m$ for any integer $n\geq 1$. 
Since all the quantities are non negative, the series $\sum_{i=0}^{\infty} E^i d$ is convergent line by line.
This implies that $\lim_{n\rightarrow\infty} E^n d=0$. If $d_{-}=\min_{1\leq i\leq p}d_i>0$, we deduce that $E^n\rightarrow 0$, element-wise.
This entails $\rho(E)<1$.

\item
  Let $\delta>0$ such that $(1+\delta)\sum_{j=1}^q(A_j+B_j)$ has a spectral radius less than one.
	If $\Gamma_{\delta}:=(1+\delta)\sum_{j=1}^q\begin{pmatrix} B_j&A_j\\B_j&A_j\end{pmatrix}$, then we also have $\rho\left(\Gamma_{\delta}\right)<1$.
	Next, from Lemma \ref{lem::poism} given in the Appendix, there exists $b>0$ such that 
 $$
\|N^{(t)}_{ f(x_1, \ldots, x_q)}\|_{r, vec} \preccurlyeq (1+\delta)\vert f(x_1, \ldots, x_q)\vert_{vec} + b\mathds{1}, 
$$
where $\mathds{1}$ denotes the vector of $\R^p$ for which all the coordinates are equal to $1$.
But, we also have 
$$
\left\vert f(x_1, \ldots, x_q)\right\vert_{\mathrm{vec}} \preccurlyeq d +\sum_{j=1}^q B_j|y_j|_{vec} + \sum_{j=1}^q A_j|s_j|_{vec} \preccurlyeq d +\sum_{j=1}^q (B_j ~ A_j)|x_j|_{vec}.
$$
We then get 

\begin{eqnarray*}
\|F(x_1, \ldots, x_q, N^{(1)})\|_{r, vec} & \preccurlyeq  & 
\begin{pmatrix}
(1+\delta)d + b\mathds{1} \\
d 
 \end{pmatrix}  + \sum_{j = 1}^q 
 \begin{pmatrix}
(1 + \delta) B_j & (1 + \delta) A_j \\
 B_j &  A_j 
 \end{pmatrix} |x_j|_{vec} \\
 & \preccurlyeq  & 
 \begin{pmatrix}
(1+\delta)d + b\mathds{1} \\
d
 \end{pmatrix} 
 + (1+\delta)\sum_{j = 1}^q 
 \begin{pmatrix}
B_j & A_j \\
 B_j &  A_j 
 \end{pmatrix} |x_j|_{vec}.\\
\end{eqnarray*}

Since $\rho\left(\Gamma_{\delta}\right)<1$, Theorem \ref{th::autoreg} (3.) applied when $\phi$ is the identity leads to the result.$\square$

\item
Assume first that $\gamma := \sum_{i=1}^q\left(\vert A_i\vert_1+\vert B_i\vert_1\right)<1$. Let $\delta>0$ to be chosen later
and $\phi(y,s)=\left(\exp\left(\delta \vert y\vert_1\right),\exp\left(\delta \vert s\vert_1\right)\right)$ for $y,s\in \R^p$. From convexity of the exponential function and  matrix norm inequalities, we can write if $t\geq q+1$:
  \begin{eqnarray*}
 \E\left[\exp(\delta \vert f(x_1,\ldots,x_q)\vert_1)\right] & \leq & \exp\left(\delta \vert d\vert_1 + \sum_{j = 1}^q \delta \vert B_j\vert_1\vert y_j\vert_1 + \delta \vert A_j\vert_1\vert s_j\vert_1 \right)\\
 &\leq & c +  \sum_{j = 1}^q  \vert B_j\vert_1 \exp(\delta\vert y_j\vert_1) +  \vert A_j \vert_1 \exp(\delta\vert s_j\vert_1), \\ c = (1-\gamma)\exp\left(\frac{\delta\vert d\vert_1}{1-\gamma}\right).\\
 \end{eqnarray*}
  Furthermore, from H\"{o}lder inequality, we have setting $\overline{\lambda}=f(x_1,\ldots,x_q)$,
  \begin{eqnarray*}
 \E\left[\exp(\delta\vert N^{(1)}_{\overline{\lambda}}\vert_1)\right] & = & \E\left[\prod_{j = 1}^p \exp\left(\delta N^{(1)}_{j,\overline{\lambda}_j} \right)\right] \\
 &\leq & \prod_{j = 1}^p \E^{1/p} \left[\exp\left(p \delta N^{(1)}_{j,\overline{\lambda}_j}\right)\right] \\
 & = &  \prod_{j = 1}^p  \left(\exp\left(\overline{\lambda}_j[\exp^{p\delta} - 1]\right)\right)^{1/p} \\
 & = & \exp\left(\vert\overline{\lambda}\vert_1 \frac{\exp^{p\delta} - 1}{p}\right).
 \end{eqnarray*}
 But $\frac{\exp^x - 1}{x} \downarrow 1$ as $x > 0$ tends to $0$. If $\epsilon > 0$ is such that $(1+\epsilon)\gamma < 1$, let us choose $\delta = \delta(\epsilon)$ such that $\exp^{p\delta} -1 \leq (1+\epsilon)\delta p$. For the couple $(\epsilon, \delta  = \delta(\epsilon))$, one can write :
\begin{eqnarray*} 
  \E\left[\exp(\delta\vert N^{(1)}_{\overline{\lambda}}\vert_1)\right] &\leq& \exp\left(\vert \overline{\lambda}\vert_1(1+\epsilon)\delta\right)\\
	&\leq& c' +  \sum_{j = 1}^q  (1+\epsilon)\vert B_j\vert_1 \exp(\delta\vert s_j\vert_1) +  (1+\epsilon)\vert A_j\vert_1 \exp(\delta\vert y_j\vert_1),
	\end{eqnarray*}
	with
 $$c' = \left(1-(1+\epsilon)\gamma\right)\exp\left(\frac{\delta(1+\epsilon)\vert d\vert_1}{1-(1+\epsilon)\gamma}\right).$$
The second inequality follows from the convexity of the exponential function, as previously. 
We then obtain
 \begin{eqnarray*}
\E\left[\phi\left(F(x,N^{(1)})\right)\right]
 & \preccurlyeq  & 
\begin{pmatrix}
c' \\
c 
 \end{pmatrix}  + \sum_{j = 1}^q 
 \begin{pmatrix}
(1 + \epsilon) \vert B_j\vert_1 & (1 + \epsilon) \vert A_j\vert_1 \\
 \vert B_j\vert_1 &  \vert A_j\vert_1 
 \end{pmatrix} 
\phi(x_j)\\
 & \preccurlyeq  & 
 \begin{pmatrix}
c' \\
c 
 \end{pmatrix}  + (1 + \epsilon) \sum_{j = 1}^q 
 \begin{pmatrix}
 \vert B_j\vert_1 &  \vert A_j\vert_1 \\
 \vert B_j\vert_1 &  \vert A_j\vert_1 
 \end{pmatrix}\phi(x_j). \\
\end{eqnarray*}
Since the eigenvalues of $\sum_{j=1}^q\begin{pmatrix}
 \vert B_j\vert_1 &  \vert A_j\vert_1 \\
 \vert B_j\vert_1 &  \vert A_j\vert_1 
 \end{pmatrix}$ are $\gamma$ et $0$, the spectral radius of  $$\Gamma_\epsilon = (1 + \epsilon) \sum_{j = 1}^q 
 \begin{pmatrix}
 \vert B_j\vert_1 &  \vert A_j\vert_1 \\
 \vert B_j\vert_1 &  \vert A_j\vert_1 
 \end{pmatrix}$$ is $(1+\epsilon)\gamma$ and is less than 1.
The result then follows from Theorem \ref{th::autoreg} (3.) applied with $r=1$.

Next, we assume that $\gamma:=\left\vert \sum_{i=1}^q\left(A_i+B_i\right)\right\vert_{\infty}<1$. Some arguments previously used yield $\forall k = 1, \ldots, p$,
 $$
  \E\left[\exp(\delta \overline{\lambda}_k)\right] \leq   e_k + \sum_{j = 1}^q \sum_{l = 1}^p B_j(k,l) \exp(\delta  y_{\ell,j}) + A_j(k,l) \exp(\delta s_{\ell,j}),
 $$ and 
 $$
  \E\left[\exp(\delta N^{(1)}_{k,\overline{\lambda}_k})\right] \leq \exp(\delta (1+\epsilon) \overline{\lambda}_{k,t}) \leq e_k' + (1+\epsilon)\sum_{j = 1}^q \sum_{l = 1}^p A_j(k,l) \exp(\delta  s_{\ell,j}) + B_j(k,l) \exp(\delta  y_{\ell,j}),
 $$
 with $e_k = \left(1-\gamma_k\right)\exp\left(\frac{\delta d_k}{1-\gamma_k}\right), ~ e_k' = \left(1-(1+\epsilon)\gamma_k\right)\exp\left(\frac{\delta(1+\epsilon)d_k}{1-(1+\epsilon)\gamma_k}\right)$ and $\gamma_k = \sum_{j = 1}^q \sum_{l = 1}^p A_j(k,l) + B_j(k,l)$ where $(\epsilon, \delta= \delta(\epsilon))$ satisfy $(1+\epsilon)\sup_k \gamma_k=(1+\epsilon)\gamma<1$ and $\exp^{\delta} -1 \leq (1+\epsilon)\delta .$ Therefore, 
 
\begin{eqnarray*}
\begin{pmatrix}
 \E\left[\exp(\delta N^{(1)}_{\overline{\gamma}} \right]_{vec} \\
\E\left[\exp(\delta \overline{\lambda}) \right]_{vec} \\
\end{pmatrix} & \preccurlyeq  & 
\begin{pmatrix}
e'_{vec} \\
e_{vec}
 \end{pmatrix}  +  (1+\epsilon) \displaystyle\sum_{j = 1}^q 
\begin{pmatrix}
B_j & A_j \\
B_j & A_j
\end{pmatrix}
\begin{pmatrix}
 \exp(\delta y_{j})_{vec}  \\
\exp(\delta s_{j})_{vec}  \\
\end{pmatrix} \\
\end{eqnarray*}
 With $M=\sum_{j=1}^q 
(A_j + B_j)$, condition  $\vert M\vert_{\infty}=\max_{1 \leq j \leq q}  \sum_{l = 1}^q M(j,l) < 1$ ensures that the spectral radius of the matrix $\sum_{j=1}^q \begin{pmatrix}
 B_j & A_j \\ B_j & A_j 
\end{pmatrix}$ is less than 1. Then, one can find $\epsilon$ such that the spectral radius of 
$$\Gamma_\epsilon' = (1+\epsilon)\sum_{j=1}^q \begin{pmatrix}
 B_j & A_j \\ B_j & A_j 
\end{pmatrix}$$
is less than 1. The result then follows from Theorem \ref{th::autoreg} (3.), setting $r=1$ and for $y,s\in\R^p$,
$$\phi((y,s))=\left(\exp(\delta \vert y_1\vert),\ldots,\exp(\delta \vert y_p\vert),\exp(\delta \vert s_1\vert),\ldots,\exp(\delta\vert s_p\vert)\right)'.\square$$

\end{enumerate}

\subsection{Log linear model}

We now consider a second model called log-linear in the literature. See in particular \citet{Fokianos2011b} for the univariate case and \citet{DFT} for the multivariate case. 
In the multivariate case, the model is defined similarly to (\ref{multlin}) except that
\begin{equation}\label{loglin} 
\lambda_t=\exp(\mu_t),\quad \mu_t=d+\sum_{j=1}^q A_j\mu_{t-j}+\sum_{j=1}^q B_j \log(\mathds{1}+Y_{t-j}).
\end{equation}
Here, the functions $\exp$ and $\log$ are applied component-wise and the matrices $A_j,B_j$ can now have negative elements.
As before, we adopt the convention of column vectors.

\begin{theo}
\label{th::llm}
Consider the $\log$ linear model (\ref{loglin}).
\begin{enumerate}
\item
 Assume that $\rho\left(\sum_{i=1}^q (\vert A_i\vert_{vec} + \vert B_i\vert_{vec})\right)< 1$. Then there exists a unique non anticipative, stationary and integrable process $(Y_t)_{t\in\Z}$ with $\E(\vert Y_t\vert_1) < \infty$.
\item
Assume that $\left\vert \sum_{i=1}^q\left(\left\vert A_i\right\vert_{vec}+\left\vert B_i\right\vert_{vec}\right)\right\vert_{\infty}<1$.
There then exists $\delta>0$ such that $\E\left(\exp\left(\delta\vert Y_0\vert_1\right)\right)<\infty$ and $\E\left(\exp\left(\delta\vert \lambda_0\vert_1\right)\right)<\infty$.
\end{enumerate}
\end{theo}

\paragraph{Notes}

\begin{enumerate}
\item
To compare our results with that of \citet{DFT}, we assume $q=1$. Using contraction properties of autoregressive processes, \citet{DFT} used the condition $\vert A_1\vert_1+\vert B_1\vert_1$ for studying existence of a stationary solution for (\ref{loglin}). Our condition in point $1$ of Theorem \ref{loglin} is weaker.
On the other hand, using perturbation methods, \citet{DFT} showed that one can approximate the stationary solution from an ergodic Markov chain when $\vert A_1\vert_2+\vert B_1\vert_2<1$. When $p=1$, the latter condition coincides with our but when $p>1$ they cannot be compared. Our condition is only guaranteed to be weaker
when the matrices $A_1$ and $B_1$ have nonnegative coefficients. We then provide a different result which complements the existing ones.

\item
We obtain directly existence of exponential moments for the solution using a contraction condition for the norm $\vert \cdot\vert_{\infty}$. 
We are not aware of such result even in the univariate case. \citet{DFT} only studied existence of polynomial moments but once again their assumptions 
are based on the norm $\vert\cdot\vert_1$ and $\vert \cdot\vert_2$ and cannot be compared directly with our, even for nonnegative matrices.
Nevertheless, for $p=1$, all the conditions are equivalent to $\vert A_1\vert+\vert B_1\vert<1$, and we improve existing results by showing existence of exponential moments.  
\end{enumerate}

\paragraph{Proof of Theorem \ref{th::llm}}
\begin{enumerate}
\item
For any solution, setting $X_t = (\log(\mathds{1} + Y_t)', \log(\lambda_t)')'$. We can write:
 $X_t = F(X_{t-1}, \ldots, X_{t-q}, N^{(t)})$ with 
$$F(x_1,\ldots,x_q)= \left(\log(\mathds{1} + N^{(t)}_{f(x_1, \ldots, x_q)})',\log(f(x_1, \ldots, x_q))' \right)',$$
 where $\forall j = 1, \ldots, q, x_j = (u_j, v_j)\in (\R_+^p \times \R^p)$ and 
$$\log(f(x_1, \ldots, x_q)) = b_0+\sum_{j=1}^q B_j u_j+\sum_{j=1}^q A_j v_j.$$
  
     For $x \in (\R_+^p \times \R^p)^q$, 
      $$
      \E\left[\left\vert F\left(x_1,\ldots,x_q,N^{(1)}\right)\right\vert_1\right] =  \mathds{1}'\left(\left\vert \log(f(x_1, \ldots, x_q))\right\vert_{vec} + \log(\mathds{1} +  f(x_1, \ldots, x_q))\right) < \infty.
      $$

Using Jensen's inequality to the function $x \mapsto \log(1+x)$ and Poisson process properties, we obtain that, for a given Poisson process $N$ (for more details, see \citet{Fokianos2011b}, proof of Lemma $2.1$), 
 $$
 \E\left(\log\left(\frac{1+N_t}{1+N_s}\right)\right) \leq \log(t)-\log(s).
 $$
 For $x = (x_1, \ldots, x_q), x' = (x_1', \ldots, x_q') \in E^q$, 
 \begin{eqnarray*}
 \E\left[\left\vert F(x,N^{(1)})-F(x',N^{(1)})\right\vert_{vec}\right] & = & 
 \begin{pmatrix}
  \E\left(\left\vert\log\left(\frac{1+N^{(1)}_{f(x)}}{1+N^{(1)}_{f(x')}}\right)\right\vert_{vec}\right) \\
  \left\vert \log(f(x)) - \log(f(x')) \right\vert_{vec} 
 \end{pmatrix} \\
 & \preceq & \begin{pmatrix}
 \left\vert \log(f(x)) - \log(f(x')) \right\vert_{vec} \\
 \left\vert \log(f(x)) - \log(f(x')) \right\vert_{vec}
 \end{pmatrix} \\
 & \preceq & 
 \sum_{j=1}^q 
 \begin{pmatrix}
\vert B_j\vert_{vec}  & \vert A_j\vert_{vec}\\
\vert B_j\vert_{vec} & \vert A_j\vert_{vec}
 \end{pmatrix} \left\vert x_j-x'_j\right\vert_{vec}, 
 \end{eqnarray*}
 where for $x\in \R^d$ and $y\in \R_*^d$, $\frac{x}{y} = (\frac{x_1}{y_1}, \ldots, \frac{x_d}{y_d})'.$ Note that the matrices  $ \Gamma = \sum_{j=1}^q \begin{pmatrix}
\vert B_j\vert_{vec} & \vert A_j\vert_{vec} \\
\vert B_j\vert_{vec} & \vert A_j\vert_{vec} 
 \end{pmatrix}$ and $\sum_{j=1}^q 
\left(\vert A_j\vert_{vec} + \vert B_j\vert_{vec}\right)
$ have the same spectral radius. The result then follows from theorem \ref{th::autoreg} (1.).

\item
 We will use Theorem \ref{th::autoreg} (3.) with $r=1$ and for $(u,v)\in \R_+^p\times \R^p$,
$$\phi(u,v)=\left(\exp\left(\delta \vert \exp(u_1)-1\vert\right),\ldots,\exp\left(\delta \vert \exp(u_p)-1\vert\right),\exp\left(\delta \exp(\vert v_1\vert)\right),\ldots,\exp\left(\delta \exp(\vert v_p\vert)\right)\right),$$
with $\delta>0$ to be specified latter.
Setting for $x_i=(u_i,v_i)\in \R_+^p\times \R^p$ for $1\leq i\leq p$, $\overline{\lambda}=f\left(x_1,\ldots,x_q\right)$ and $\overline{\mu}=\log\left(\overline{\lambda}\right)$.
 We have for $1\leq k\leq p$,
$$\exp(|\overline{\mu}_k|) \leq   e_k + \sum_{j = 1}^q \sum_{l=1}^p |B_j(k,l)|\left[\exp\left(u_{\ell,j}\right)-1\right] + |A_j(k,l)|\exp(\vert v_{\ell,j}\vert),$$
 and 
 $$\E\left[\exp\left(\delta \exp(\vert \overline{\mu}_k\vert)\right)\right] \leq   e_k' + \sum_{j=1}^q \sum_{l=1}^p |B_j(k,l)|\exp\left(\delta \left[\exp(u_{\ell,j})-1\right]\right)+|A_j(k,l)|\exp\left(\delta \exp(\vert v_{\ell,j}\vert)\right),$$ 

with $$e_k =\left(1-\gamma_k\right)\exp\left(\frac{ d_k}{1-\gamma_k}\right) + \sum_{j = 1}^q \sum_{l = 1}^p |B_j(k,l)|, ~ e_k' = \left(1-\gamma_k\right)\exp\left(\frac{\delta e_k}{1-\gamma_k}\right)$$ and $\gamma_k = \sum_{j = 1}^q \sum_{l = 1}^p |A_j(k,l)| + |B_j(k,l)|$.  We also have
\begin{eqnarray*}
\E\left[\exp\left(\delta N^{(1)}_{k,\overline{\lambda}_k}\right)\right] &\leq &\exp(\delta (1+\epsilon) \overline{\lambda}_k)\\
& \leq& \tilde{e}_k + (1+\epsilon)\sum_{j = 1}^q \sum_{l = 1}^p |B_j(k,l)|\exp\left(\delta \left[\exp(u_{\ell,j})-1\right]\right) + |A_j(k,l)| \exp\left(\delta \exp(\vert v_{\ell,j}\vert)\right),
\end{eqnarray*}
with $\tilde{e}_k=\left(1-(1+\epsilon)\gamma_k\right)\exp\left(\frac{\delta(1+\epsilon)e_k}{1-(1+\epsilon)\gamma_k}\right)$ where $(\epsilon, \delta= \delta(\epsilon))$ satisfy $(1+\epsilon)\sup_k \gamma_k \leq 1$ and $\exp^{\delta} -1 \leq (1+\epsilon)\delta .$ Therefore, 
 
\begin{eqnarray*}
\E\left[\phi\left(F(x,N^{(1)})\right)\right]
& \preccurlyeq  & 
\begin{pmatrix}
\tilde{e}_{vec} \\
e'_{vec}
 \end{pmatrix}  +  (1+\epsilon) \displaystyle\sum_{j = 1}^q 
\begin{pmatrix}
\vert B_j\vert_{vec} & \vert A_j\vert_{vec} \\
\vert B_j\vert_{vec} & \vert A_j\vert_{vec}
\end{pmatrix}
\begin{pmatrix}
 \exp\left(\delta \left(\exp(u_j)-\mathds{1}\right)\right) \\
\exp\left(\delta \exp\left(\vert v_j\vert_{vec}\right)\right)\\
\end{pmatrix} \\
\end{eqnarray*}
 With $M=\sum_{j=1}^q 
( \vert B_j\vert_{vec} + \vert A_j\vert_{vec})$, condition  $\vert M\vert_{\infty}<1$ ensures that the spectral radius of the matrix $\sum_{j=1}^q \begin{pmatrix}
 \vert B_j\vert_{vec} & \vert A_j\vert_{vec}\\ \vert B_j\vert_{vec} & \vert A_j\vert_{vec} 
\end{pmatrix}$ is less than 1. Then, one can find $\epsilon$ such that the spectral radius of 
$$\Gamma_\epsilon = (1+\epsilon)\sum_{j=1}^q \begin{pmatrix}
 \vert B_j\vert_{vec} & \vert A_j\vert_{vec} \\ \vert B_j\vert_{vec} & \vert A_j\vert_{vec} 
\end{pmatrix}$$
is less than 1. Theorem \ref{th::autoreg} (3.) then leads to the result. $\square$
\end{enumerate}

\section{Appendix}\label{APP}
For a square matrix $G$, we denote by $\rho(G)$ its spectral radius.

\begin{lem}\label{calc}
Let $E_1,\ldots,E_q$ be square matrices of size $e\times e$, with nonnegative elements and such that $\rho\left(E_1+\cdots+E_q\right)<1$.
\begin{enumerate}
\item
If $F$ denotes the companion matrix associated to $E_1,\ldots,E_q$, i.e. 
$$F=\begin{pmatrix} E_1&E_2&\cdots& E_q\\ &I_{(q-1)e}& &0_{(q-1)e,e}\end{pmatrix},$$
then $\rho(F)<1$.
\item
Let also $(v_n)_{n\geq 1}$ a sequence of vectors of $\R_+^e$  and $b\in\R_+^e$ such that 
$$v_n\preceq \sum_{i=1}^q E_i v_{n-i}+b,\quad i\geq q+1.$$
Let $\vert\cdot\vert_1$ be the $\ell_1-$norm on $\R^e$. There exists $C,D>0$ and $\overline{\rho}\in (0,1)$, not depending on $(v_n)_{n\geq 1}$, such that $\vert v_n\vert_1\leq C\overline{\rho}^n\sum_{i=1}^q\vert v_i\vert_1+\frac{C \vert b\vert_1}{1-\overline{\rho}}$.
\end{enumerate}
\end{lem}

\paragraph{Proof of Lemma \ref{calc}}
\begin{enumerate}
\item
Let $E=E_1+\cdots+E_q$. Since $\rho(E)<1$, we have $E^n\rightarrow 0$. Suppose that $\lambda$ is an eigenvalue of $F$ of modulus greater than $1$. If $v=(v_1',\ldots,v_q')'\in \R^{qe}\setminus\{0\}$ is such that 
$F v=\lambda v$, we have $\lambda v_1=\sum_{i=1}^q E_iv_i=\sum_{i=1}^q \lambda^{1-i}E_i v_1$. This yields to 
$$\vert v_1\vert_{vec}\preceq\sum_{i=1}^q \vert\lambda\vert^{-i} E_i\vert v_1\vert_{vec}\preceq E\vert v_1\vert_{vec}\preceq E^n\vert v_1\vert_{vec}.$$
Letting $n\rightarrow \infty$, we get $v_1=0$ and then $v_2,\ldots,v_q=0$. This contradicts $v\neq 0$. Hence $\vert \lambda\vert<1$ and then $\rho(F)<1$.

\item
Setting for $n\geq q$, $u_n=\left(v_n',\ldots,v_{n-q+1}'\right)'$ and $B=\left(b',0_{1,(q-1)e}\right)'$, we have 
$$u_n\preceq F u_{n-1}+B\preceq F^{n-q}u_q+\sum_{i=0}^{n-q-1}F^iB,$$ with $F$ being the companion matrix 
associated to the matrices $E_1,\ldots,E_q$ and which is defined in the previous point.
We still denote by $\vert\cdot\vert_1$ the $\ell_1-$norm on $\R^{eq}$.
From the previous point, we have $\rho(F)<1$ and then if $\varepsilon>0$ is such that $\overline{\rho}=\rho(F)+\varepsilon<1$, we have $\left\vert F^n\right\vert_1\leq C\overline{\rho}^n$ for $C>0$ only depending on the matrix $F$. 
Then if $n\geq q+1$, 
$$\vert v_n\vert_1\leq \vert u_n\vert_1\leq C\overline{\rho}^{n-q}\vert u_q\vert_1+C\sum_{i=0}^{n-q-1}\overline{\rho}^i\vert B\vert_1.$$
Since $\vert B\vert_1=\vert b\vert_1$ and $\vert u_q\vert_1\leq \sum_{i=1}^q \vert v_i\vert_1$, this leads to the result.$\square$
\end{enumerate}

\begin{lem}
\label{lem::poism}
Let $\lambda>0$ and $X_{\lambda}$ $\mathrm{Poisson}$ variable  with parameter $\lambda$.   Then, $\forall r \geq 1$ and any $\delta\in (0,1)$,  there exists $b_{r,\delta}$, not depending on $\lambda$ and such that
$$
\|X_{\lambda}\|_r \leq (1+\delta)\lambda+ b_{r,\delta}.
$$
\end{lem}

\paragraph{Proof of lemma \ref{lem::poism}}
We have the equality $\mathbb{E}(X_{\lambda}^r) = \sum_{i = 1}^r \lambda^i \left\{\begin{array}{c}
   r  \\
    i 
\end{array} \right\}  $ with  $\left\{\begin{array}{c}
    r  \\
    i 
\end{array} \right\}$ are the  Sterling's numbers of second kind. 
See for instance \cite{johnson2005univariate}.

Then 
$$
\mathbb{E}(X_{\lambda}^r)  = \lambda^r + \sum_{i = 1}^{r-1} \lambda^i \left\{\begin{array}{c}
    r  \\
    i 
\end{array}\right\} 
\leq  \lambda^r + C_r(\lambda + \lambda^{r-1}),
$$
where $C_r>0$ only depends on $r$.
But, we can notice that,  for any $\delta>0$, there exists $\exists \overline{b}_{\delta, r}>0$ such that for all $x\geq 0$ : $x + x^{r-1} \leq \delta' x^r + \overline{b}_{\delta, r}$ with $\delta'=\frac{(1+\delta)^r-1}{C_r}$. Then $\mathbb{E}(X^r) \leq (1+C_r\delta')\lambda^r + C_r \overline{b}_{\delta, r}$. 
Therefore $\|X\|_r \leq (1+C_r\delta')^{1/r} \lambda + C_r^{1/r}\overline{b}_{\delta, r}^{1/r}$.
Setting $b_{\delta,r}=C_r^{1/r}\overline{b}_{\delta, r}^{1/r}$, we get the result.$\square$

\bibliographystyle{plainnat}
\bibliography{biblio}

\end{document}